\documentclass[11pt]{amsart}

\usepackage{etex}
\usepackage{amsmath, amssymb}
\usepackage{array}
\usepackage[frame,cmtip,arrow,matrix,line,graph,curve]{xy}
\usepackage{graphpap, color, paralist, pstricks}
\usepackage[mathscr]{eucal}
\usepackage[pdftex]{graphicx}
\usepackage[pdftex,colorlinks,backref=page,citecolor=blue]{hyperref}
\usepackage{tikz}

\newtheorem{theorem}{Theorem}[section]

\newtheorem{lemma}{Lemma}[section]
\newtheorem{proposition}{Proposition}[section]
\newtheorem{definition}{Definition}[section]

\newtheorem{remark}{Remark}[section]

\newcommand{\cc}{\mathbb{C}}
\newcommand{\pp}{\mathbb{P}}
\newcommand{\zz}{\mathbb{Z}}

\newcommand{\bM}{\mathbf{M} }
\newcommand{\bT}{\mathbf{T} }

\newcommand{\id}{\textrm{id}}

\newcommand{\EEnd}{\mathscr E\!nd}
\newcommand{\HHom}{\mathscr H\!om}

\newcommand{\pr}{\mathrm{pr}}

\newcommand{\Jac}{\mathrm{Jac}}
\newcommand{\im}{\mathrm{im \,}}
\newcommand{\coker}{\mathrm{coker \,}}

\newcommand{\Spec}{\mathrm{Spec} \,}
\newcommand{\sym}{\mathrm{Sym}}
\newcommand{\PMod}{\mathrm{PMod}}
\newcommand{\Nm}{\mathrm{Nm}}
\newcommand{\Qu}{\mathrm{Quot}}
\newcommand{\rk}{\mathrm{rank} \,}
\newcommand{\supp}{\mathrm{Supp}\,}

\newcommand{\cC}{\mathcal{C} }
\newcommand{\cD}{\mathcal{D} }
\newcommand{\cE}{\mathcal{E} }
\newcommand{\cF}{\mathcal{F} }
\newcommand{\cH}{\mathcal{H} }
\newcommand{\cL}{\mathcal{L} }
\newcommand{\cM}{\mathcal{M} }
\newcommand{\cO}{\mathcal{O} }
\newcommand{\cP}{\mathcal{P} }
\newcommand{\cR}{\mathcal{R} }
\newcommand{\cT}{\mathcal{T} }
\newcommand{\cU}{\mathcal{U} }

\newcommand{\tT}{\tilde{T} }
\newcommand{\tX}{\tilde{X} }
\newcommand{\tf}{\tilde{f} }

\newcommand{\tlt}{\tilde{t} }

\newcommand{\ty}{\tilde{y} }
\newcommand{\tp}{\tilde{p} }
\newcommand{\tv}{\tilde{v} }
\newcommand{\tmu}{\tilde{\mu}}
\newcommand{\tpi}{\tilde{\pi}}
\newcommand{\tcC}{\tilde{\mathcal{C}} }
\newcommand{\tcL}{\tilde{\mathcal{L}} }

\newcommand{\ve}{\varepsilon }
\newcommand{\s}{\sigma }
\newcommand{\tbe}{\tilde{\beta} }

\newcommand{\ovf}{\overline{f}}

\def\Pic{\mathrm{Pic}}

\begin{document}

\title[Hecke cycles associated to rank 2 twisted Higgs bundles on a curve]{Hecke cycles associated to rank 2 twisted Higgs bundles on a curve}
\date{\today}

\author{Sang-Bum Yoo}
\address{\parbox{\linewidth}{Department of Mathematics Education, Gongju National University of Education,\\ Gongju-si, Chungcheongnam-do, 32553, Republic of Korea}}
\email{sbyoo@gjue.ac.kr}

\begin{abstract}
Let $X$ be a smooth complex projective curve of genus $g$ and let $L$ be a line bundle on $X$ with $\deg L>0$. Let $\bM$ be the moduli space of semistable rank 2 $L$-twisted Higgs bundles with trivial determinant on $X$. Let $\bM_{X}$ be the moduli space of stable rank 2 $L$-twisted Higgs bundles with determinant $\cO(-x)$ for some $x\in X$ on $X$. We construct a cycle in the product of a stack of rational maps from nonsingular curves to $\bM_{X}$ and $\Pic^{2\deg L}(X)$ by using Hecke modifications of a stable $L$-twisted Higgs bundle in $\bM$.
\end{abstract}


\maketitle

\section{Introduction}

In modern algebraic geometry, the comparison of various birational moduli spaces is a significant problem. One way to study this problem is to construct a sequence of blow-ups and blow-downs from one to another. Several successful studies in this direction can be found in \cite{CCK,CK,KL,KM}.

Our work in this paper is motivated by \cite{CCK} and \cite{NR2}. The authors of \cite{CCK} dealt with a comparison between the moduli space $\cM$ of rank 2 semistable vector bundles with trivial determinant on a smooth complex projective curve $X$ of genus $g\ge3$ and the moduli space of Hecke cycles $\cH$. $\cH$ was originally constructed in \cite{NR2} as a Zariski closure of the locus of good Hecke cycles in a Hilbert scheme of $\cM_{X}$, where $\cM_{\cO_{X}(-x)}$ is the moduli space of rank 2 stable vector bundles with determinant $\cO_{X}(-x)$ on $X$ and $\displaystyle\cM_{X}=\bigcup_{x\in X}\cM_{\cO_{X}(-x)}$. A good Hecke cycle is given by a collection of Hecke modifications of a stable vector bundle in $\cM$. The authors of \cite{NR2} proved that the stable locus of $\cM$ is isomorphic to the open dense locus of $\cH$ consisting of good Hecke cycles. Moreover they showed that $\cH$ is nonsingular. The authors of \cite{CCK} showed that $\cH$ can be obtained by a sequence of three blow-ups and a single blow-down starting from $\cM$.

\subsection{Main result}

Let $X$ be a smooth complex projective curve of genus $g$. Let $L$ be a line bundle on $X$ with $\deg L>0$. Let $\bM$ be the moduli space of semistable rank 2 $L$-twisted Higgs bundles with trivial determinant on $X$ and let $\displaystyle\bM_{X}=\bigcup_{x\in X}\bM_{\cO_{X}(-x),L}$, where $\bM_{\cO_{X}(-x),L}$ denotes the moduli space of rank 2 stable $L$-twisted Higgs bundles with determinant $\cO_{X}(-x)$ on $X$.

The purpose of this paper is to construct a cycle in the product of a stack of rational maps from nonsingular curves to $\bM_{X}$ and $\Pic^{2\deg L}(X)$. Here this cycle, which is called a Hecke cycle associated to a stable $L$-twisted Higgs bundles, is given by a collection of Hecke modifications of a stable $L$-twisted Higgs bundle in $\bM$.

For a fixed point $x\in X$, Proposition \ref{holomorphicity of phi'} and Lemma \ref{stability of elm trans} tell us that a stable $L$-twisted Higgs bundle $(E,\phi)$ in $\bM$ and $v\in\pp(E_{x}^{\vee})$ give another stable $L$-twisted Higgs bundle $(E',\phi')$ in $\bM_X$, which is called the {\bf Hecke modification of $(E,\phi)$ determined by $v$ at $x$}. It follows from Proposition \ref{a correspondence between HM on X and HM on the normalization of the spectral cover} that a Hecke modification of $(E,\phi)$ at $x$ corresponds to a point of zero dimensional subscheme $P_{(E,\phi),x}$ of $\pp(E_{x}^{\vee})=\pp^1$ of length $2$ under the assumption that the spectral curve $X_{\det\phi}$ is irreducible, the rank $1$ torsion free sheaf $\cL_{(E,\phi)}$ on $X_{\det\phi}$ associated to $(E,\phi)$ is locally free and $x\in X$ is not a branch point of the spectral cover $X_{\det\phi}\to X$. Here $P_{(E,\phi),x}$ is the projectivization of the fiber of $\cL_{(E,\phi)}^{\vee}$ over such $x\in X$. Then we get a morphism $\theta_{(E,\phi),x}:P_{(E,\phi),x}\to\bM_{\cO(-x),L}$ given by $v\mapsto (E',\phi')$. Varing $x$ in the complement of the locus of branch points of the spectral cover $X_{\det\phi}\to X$ in $X$, we get a rational map $\theta_{(E,\phi)}:P_{(E,\phi)}\dashrightarrow\bM_{X}$ over $X$, where $P_{(E,\phi)}$ is a nonsingular curve over $X$ with the fiber $P_{(E,\phi),x}$ over $x\in X$. The following is the first main result of this paper.

\begin{theorem}[Theorem \ref{Main theorem1:map theta}]\label{Main result1}
$\theta_{(E,\phi)}:P_{(E,\phi)}\dashrightarrow\bM_{X}$ is a closed immersion on the complement of the locus of branch points of the spectral cover $X_{\det\phi}\to X$ in $X$.
\end{theorem}

Now we consider a stack $\cR_{\bM_X}$ of rational maps from nonsingular curves to $\bM_X$ over the category $(Sch^{sf}/\cc)$ of normal separated schemes of finite type over $\cc$. Then $\theta_{(E,\phi)}:P_{(E,\phi)}\dashrightarrow\bM_{X}$ can be considered as an object in $\cR_{\bM_X}$. Since these rational maps $\theta_{(E,\phi)}$ are functorial, we can define a morphism of stacks $\Phi$ from an open subscheme $\bM'$ of $\bM^{s}$ into $\cR_{\bM_X}\times\Pic^{2\deg L}(X)$ given by $(E,\phi)\mapsto(\theta_{(E,\phi)},D_{\det\phi})$, where $D_{\det\phi}$ denotes the divisor of $\det\phi$. The following is the second main result of this paper.

\begin{theorem}[Theorem \ref{Main theorem2:map Phi}]\label{Main result2}
$\Phi$ is injective on objects up to isomorphism.
\end{theorem}

Here the object $\Phi((E,\phi))$ as a geometric point of the subcategory of images of $\Phi$ is called a \textbf{Hecke cycle associated to $(E,\phi)$}.

\subsection{Outline of the proofs of Theorem \ref{Main result1} and Theorem \ref{Main result2}}

For the proof of Theorem \ref{Main result1}, we show that $\theta_{(E,\phi)}$ and its differential are injective on the complement of the locus of branch points of the spectral cover $X_{\det\phi}\to X$ in $X$.

For the proof of Theorem \ref{Main result2}, we show that $\Phi(T)$ is injective on objects up to isomorphism for each scheme $T\in(Sch^{sf}/\cc)$, that is, for every pair of $(\cE,\varphi)$ and $(\cF,\psi)$ in $\bM'(T)$, if $\Phi(T)((\cE,\varphi))\cong\Phi(T)((\cF,\psi))$, then there exists an isomorphism $\ovf:X\times T\to X\times T$ and a line bundle $\cL$ on $T$ such that $(\cE,\varphi)\cong\ovf^{*}(\cF,\varphi)\otimes \pr_{T}^{*}\cL$, where $\pr_{T}:X\times T\to T$ is the projection onto $T$.

\subsection{Structure of the paper}

This paper is organized as follows. In Section \ref{ModSpHiggsBdls}, we review some of basics on twisted Higgs bundles, the spectral data associated to them and parabolic modules. In Section \ref{HeckeModHiggs}, we define a Hecke modification of a twisted Higgs bundle on $X$ and then investigate a correspondence between a Hecke modification of a twisted Higgs bundle $(E,\phi)$ on $X$ and that of $\cL_{(E,\phi)}$ on $X_{\det\phi}$. In Section \ref{Cycles}, we define the rational map $\theta_{(E,\phi)}:P_{(E,\phi)}\dashrightarrow\bM_{X}$ and the morphism of stacks $\Phi:\bM'\to\cR_{\bM_X}\times\Pic^{2\deg L}(X)$. Then we prove Theorem \ref{Main result1} and Theorem \ref{Main result2}.

\subsection*{Notations}

Throughout this paper, $X$ denotes a smooth complex projective curve of genus $g$ and $L$ denotes a line bundle on $X$ with $\deg L>0$.

\section{Preliminaries}\label{ModSpHiggsBdls}

In this section, we review some of basics on twisted Higgs bundles on $X$, the spectral data associated to them and parabolic modules. For details of the results in this section, see \cite{BNR,C93,C9394,GO,Hit,HP,Nit,Sim}.

\subsection{Twisted Higgs bundles}\label{twisted Higgs bundles}

An $L$-\textbf{twisted Higgs bundle} on $X$ is a pair of a rank $2$
vector bundle $E$ on $X$ and a section $\phi$ of
$H^0(X,\EEnd(E)\otimes L)$. For a fixed line bundle $\Lambda$ on $X$, an $L$-\textbf{twisted Higgs bundle with determinant} $\Lambda$ on $X$ is a pair of a rank $2$
vector bundle $E$ with determinant $\Lambda$ on $X$ and a section $\phi$ of
$H^0(X,\EEnd_0(E)\otimes L)$, where $\EEnd_0(E)$ denotes the traceless part of $\EEnd(E)$. Here $\phi$ is called a \textbf{Higgs field}. The notion of stability has to be imposed to construct a separated moduli space.

\begin{definition}[\cite{Hit,Nit,Sim}]
An $L$-twisted Higgs bundle $(E,\phi)$ is \textbf{stable} (respectively, \textbf{semistable}) if for any $\phi$-invariant line subbundle $F$ of $E$, we have
$$\deg(F)<\frac{\deg(E)}{2}\quad(\text{respectively, }\le).$$
\end{definition}

Two semistable $L$-twisted Higgs bundles are said to be \textbf{S-equivalent} if they have same factors in their Jordan-H\"{o}lder filtrations.

Let $\Lambda$ be a line bundle on $X$. Let $\bM_{\Lambda,L}$ be the coarse moduli space of S-equivalence classes of semistable $L$-twisted Higgs bundles with determinant $\Lambda$ on $X$. $\bM_{\Lambda,L}^{s}$ denotes the open subvariety of $\bM_{\Lambda,L}$ parametrizing stable $L$-twisted Higgs bundles. When $\Lambda=\cO_{X}$, $\bM_{\Lambda,L}$ (respectively, $\bM_{\Lambda,L}^{s}$) is denoted by $\bM$ (respectively, $\bM^{s}$). N. Nitsure proved the following result.

\begin{theorem}[\cite{Nit}]
\begin{enumerate}
\item $\bM_{\Lambda,L}$ is a connected quasi-projective variety.

\item If $\deg L\ge 2g-2$, then $\bM_{\Lambda,L}^{s}$ is a smooth irreducible variety of dimension $3\deg L$.
\end{enumerate}
\end{theorem}

\subsection{The spectral curve and BNR correspondence}

There exists a relation between spectral curves and $L$-twisted Higgs bundles. See \cite{BNR,HP,Sim} for details.

Let $\bT=\mathbf{Spec}(\sym(L^{-1}))$ and let $\pi:\bT\to X$ be the projection. We denote the tautological section of $\pi^{*}L$ by $\lambda$. For $s\in H^0(X,L^{2})$, the \textbf{spectral curve} $X_{s}$ associated to $s$ is the zero scheme in $\bT$ of the section
$$\lambda^{2}+\pi^{*}s\in H^0(\bT,\pi^{*}L^{2}).$$
The spectral curve $X_{s}$ is reduced for nonzero section $s$ (See Remark 3.2 of \cite{GO}), but it can be singular and reducible. Moreover, the map $\pi:X_{s}\to X$ induced from $\pi:\bT\to X$ is a double cover of $X$, ramified over $D_{s}:=\mathrm{div}(s)\in\sym^{2\deg L}(X)$. The double cover $\pi:X_{s}\to X$ is usually called a \textbf{spectral cover}. We have the following simple criterions for the smoothness and irreducibility of $X_{s}$.

\begin{remark}[Example 3.4 and Remark 3.5 of \cite{BNR}]\label{smooth and irreducible}
\begin{enumerate}
\item $X_{s}$ is smooth if and only if $s$ has only simple zeros.

\item $X_{s}$ is irreducible if and only if $s$ is not the square of a section of $L$.
\end{enumerate}
\end{remark}

There is a correspondence between data on spectral curves and $L$-twisted Higgs bundles as follows.

\begin{theorem}[Proposition 3.6 of \cite{BNR}, Proposition 6.1 of \cite{HP} and Lemma 6.8 of \cite{Sim}]\label{BNR correspondence}
For $s\in H^0(X,L^{2})$, there is a bijective correspondence between isomorphism classes of semistable rank $1$ torsion free sheaves $\cL$ on $X_s$ and isomorphism classes of semistable $L$-twisted Higgs bundles $(E,\phi)$ with $\det\phi=s$. The correspondence is given by associating to any line bundle $\cM$ on $X_s$ the sheaf $\pi_{*}\cM$ on $X$ and the natural morphism
$$\pi_{*}\cM\to\pi_{*}\cM\otimes L\cong\pi_{*}(\cM\otimes\pi^{*}L)$$
given by the tautological section $\lambda$ of $\pi^{*}L$.
\end{theorem}

\subsection{Parabolic modules on the normalization of the spectral curve}\label{Parabolic modules on the normalization of the spectral curve}

Suppose that $X_{s}$ is singular and irreducible for some $s\in H^0(X,L^{2})$. Then it follows from Remark \ref{smooth and irreducible} that $s\in H^0(X,L^{2})$ has at least one multiple zero, but it is not a square of a section of $L$.

Assume that $s\in H^0(X,L^{2})$ be a section with $r_{1}$ zeros with even multiplicities $m_{i}$ ($i=1,\cdots,r_{1}$) and $r_{2}$ zeros with odd multiplicities $m'_{j}$ ($j=1,\cdots,r_{2}$). Then $X_{s}$ has $r_{1}$ nodes $x_{1},\cdots,x_{r_1}$ of types $A_{m_{i}-1}$ ($i=1,\cdots,r_{1}$) and $r_{2}$ cusps $y_{1},\cdots,y_{r_2}$ of types $A_{m'_{j}-1}$ ($j=1,\cdots,r_{2}$).

Let $\tpi:\tX_{s}\to X_{s}$ be the normalization. Most of the rank $1$ torsion free sheaves $X_{s}$ can be identified with a line bundle on $\tX_{s}$ with an additional structure reflecting informations of gluings via the normalization $\tpi:\tX_{s}\to X_{s}$. The additional structure is defined as follows.

\begin{definition}[Definition 5.1 of \cite{GO}]
Let $\tpi^{-1}(x_{i})=\{x_{i1},x_{i2}\}$ and let $\tpi^{-1}(y_{j})=\{\ty_{j}\}$. A \textbf{parabolic module on} $\tX_{s}$ is a pair $(\cL,U(\cL))$, where
\begin{enumerate}
\item $\cL\in\Jac(\tX_{s})$;

\item $U(\cL)=(U_{1}(\cL),\cdots,U_{r_1}(\cL),U'_{1}(\cL),\cdots,U'_{r_2}(\cL))$, where
\begin{enumerate}
\item for each $i=1,\cdots,r_{1}$, $U_{i}(\cL)$ is an $m_{i}/2$-dimensional subspace of the vector space $(\cL_{x_{i1}}\oplus\cL_{x_{i2}})^{m_{i}/2}$ which is also an $\cO_{X_{s},x_{i}}$-module via $\tpi_{*}$;

\item for each $j=1,\cdots,r_{2}$, $U'_{j}(\cL)$ is an $(m'_{j}-1)/2$-dimensional subspace of the vector space $\cL_{\ty_{j}}^{m'_{j}-1}$ which is also an $\cO_{X_{s},y_{j}}$-module via $\tpi_{*}$.
\end{enumerate}
\end{enumerate}
\end{definition}

P.R. Cook constructed the moduli space $\PMod_{\vec{m}}(\tX_{s})$ of parabolic modules on $\tX_{s}$, associated to $r_{1}+r_{2}$ singularities of $X_{s}$ of type indexed by $\vec{m}=(m_{1},\cdots,m_{r_1},m'_{1},\cdots,m'_{r_2})$ in \cite{C93} and \cite{C9394}. We have a finite morphism
$$\tau:\PMod_{\vec{m}}(\tX_{s})\to\overline{\Jac(X_{s})},$$
where $\tau(\cL,U(\cL))$ is the kernel of the restriction
$$\tpi_{*}\cL\to\bigoplus_{i=1}^{r_1}(\cL_{x_{i1}}\oplus\cL_{x_{i2}})^{m_{i}/2}/U_{i}(\cL)\oplus\bigoplus_{j=1}^{r_2}\cL_{\ty_{j}}^{m'_{j}-1}/U'_{j}(\cL)$$
and $\overline{\Jac(X_{s})}$ be the compactified Jacobian of $X_{s}$ parametrizing rank $1$, degree $0$ torsion free sheaves on $X_{s}$. P.R. Cook proved the following result in \cite{C93}.

\begin{proposition}[Theorem 4.4.1 of \cite{C93}]\label{identification via tau}
The restriction $\tau|_{\tau^{-1}(\Jac(X_{s}))}$ gives an isomorphism $\tau^{-1}(\Jac(X_{s}))\cong\Jac(X_{s})$.
\end{proposition}
\begin{proof}
For the readers' convenience, we present the proof of the statement that $\tau|_{\tau^{-1}(\Jac(X_{s}))}$ gives an isomorphism $\tau^{-1}(\Jac(X_{s}))\cong\Jac(X_{s})$ pointwisely. We follow the proofs of Lemma 3.4.3 and Corollary 3.4.4 of \cite{C93}. Let $\displaystyle D=\sum_{i=1}^{r_1}\frac{m_i}{2}x_{i1}+\sum_{i=1}^{r_1}\frac{m_i}{2}x_{i2}+\sum_{j=1}^{r_2}\frac{m'_{j}-1}{2}\ty_{j}$.

We first claim that $\cM\cong\tau(\tau^{-1}(\cM))$ for any $\cM\in\Jac(X_{s})$. For $\cM\in\Jac(X_{s})$, set $\tau^{-1}(\cM)=(\tpi^{*}\cM,U(\tpi^{*}\cM))$, where $U(\tpi^{*}\cM)$ is the subbundle $\cM\otimes_{\cO_{X_{s}}}\bigg(\frac{\cO_{X_{s}}}{\tpi_{*}\cO_{\tX_{s}}(-D)}\bigg)$ of
$$\bigg(\cM\otimes_{\cO_{X_{s}}}\bigg(\frac{\cO_{X_{s}}}{\tpi_{*}\cO_{\tX_{s}}(-D)}\bigg)\bigg)\otimes_{\cO_{X_{s}}}\tpi_{*}\cO_{\tX_{s}}=\tpi^{*}\cM\otimes_{\cO_{\tX_{s}}}\frac{\cO_{\tX_{s}}}{\cO_{\tX_{s}}(-D)}.$$
Then we have the following short exact sequence
$$0\to\tau(\tau^{-1}(\cM))\to\tpi_{*}\tpi^{*}\cM\to\frac{\tpi^{*}\cM\otimes_{\cO_{\tX_{s}}}\bigg(\cO_{\tX_{s}}/\cO_{\tX_{s}}(-D)\bigg)}{\cM\otimes_{\cO_{X_{s}}}\bigg(\cO_{X_{s}}/\tpi_{*}\cO_{\tX_{s}}(-D)\bigg)}$$
$$\quad\quad\quad\quad\quad\quad\quad\quad\quad\quad\quad\quad\quad\quad\quad=\cM\otimes_{\cO_{X_{s}}}\tpi_{*}\cO_{\tX_{s}}/\cO_{X_{s}}\to0.$$
On the other hand, we have a canonical injection $\cM\to\tpi_{*}\tpi^{*}\cM$ whose cokernel is $\cM\otimes_{\cO_{X_{s}}}\tpi_{*}\cO_{\tX_{s}}/\cO_{X_{s}}$. Thus there is a morphism $\cM\to\tau(\tau^{-1}(\cM))$. Since $\deg\cM=\deg\tau(\tau^{-1}(\cM))$, this is an isomorphism.

Next we claim that $\tau^{-1}(\tau((\cL,U(\cL))))=(\cL,U(\cL))$ for any $(\cL,U(\cL))\in\tau^{-1}(\Jac(X_{s}))$. Let $\cM=\tau((\cL,U(\cL)))$. Since $\deg\tpi^{*}\cM=\deg\cL$, we have $\tpi^{*}\cM\cong\cL$. We have only to see that $U(\cL)=\cM\otimes_{\cO_{X_{s}}}\bigg(\frac{\cO_{X_{s}}}{\tpi_{*}\cO_{\tX_{s}}(-D)}\bigg)$. Since $\cM$ fits into the following short exact sequences
$$0\to\cM\to\tpi_{*}\tpi^{*}\cM\to\frac{\tpi^{*}\cM\otimes_{\cO_{\tX_{s}}}\bigg(\cO_{\tX_{s}}/\cO_{\tX_{s}}(-D)\bigg)}{U(\cL)}\to0$$
and
$$0\to\cM\to\tpi_{*}\tpi^{*}\cM\to\cM\otimes_{\cO_{X_{s}}}\tpi_{*}\cO_{\tX_{s}}/\cO_{X_{s}}\to0,$$
we get $U(\cL)=\cM\otimes_{\cO_{X_{s}}}\bigg(\frac{\cO_{X_{s}}}{\tpi_{*}\cO_{\tX_{s}}(-D)}\bigg)$.

Following the proof of Theorem 4.3.4 of \cite{C93}, we complete the proof.
\end{proof}

Consider the projection on the first factor
$$\pr_{1}:\PMod_{\vec{m}}(\tX_{s})\to\Jac(\tX_{s}),\quad(\cL,U(\cL))\mapsto\cL.$$
Under the identification of Proposition \ref{identification via tau}, $\pr_{1}|_{\tau^{-1}(\Jac(X_{s}))}$ can be identified with $\tpi^{*}:\Jac(X_{s})\to\Jac(\tX_{s})$.

\section{Hecke modification of $L$-twisted Higgs bundles}\label{HeckeModHiggs}

The Hecke modificaion of vector bundles has been studied in various contexts. If we consider the similar construction of the Hecke modification of twisted Higgs bundles, poles may appear in its Higgs field (See the section 4.5 of \cite{Wit}). Recently, G. Wilkin gave a definition of the Hecke modification of Higgs bundles and a correspondence between Hecke modifications of Higgs bundles and those of their spectral data when the associated spectral curve is smooth (See \cite{Wil}). In the first subsection, we first give a definition of the Hecke modification of twisted Higgs bundles. Then we review the correspondence of G. Wilkin and give a generalized correspondence between Hecke modifications of twisted Higgs bundles and those of their spectral data when the associated spectral curve is possibly singular and irreducible. In the second subsection, we present some results about $(k,l)$-stability and Hecke modifications of twisted Higgs bundles, which are modified versions of some results in \cite{NR1} and \cite{NR2}.
All of these will be useful in Section \ref{Cycles}.

\subsection{A correspondence between Hecke modifications of twisted Higgs bundles and those of their spectral data}

Assume that a rank $2$ vector bundle $E$, a point $x\in X$ and a nonzero vector $v\in E_{x}^{\vee}$ are given. Then we have a surjective homomorphism $E\to\cc_{x}$, denoted by $v$, given by the composition of the restriction $E\to E_{x}$ and a quotient $v:E_{x}\to\cc_{x}$ obtained from $v\in E_{x}^{\vee}$. It is easy to see that $\ker v$ is also a rank $2$ vector bundle. $\ker v$ is called the \textbf{Hecke modification of }$E$\textbf{ determined by }$v$ at $x$. Note that $\pp(E_{x}^{\vee})$ parametrizes Hecke modifications of $E$ at $x$.

\begin{definition}\label{Hecke modification}
Let $(E,\phi)$ be an $L$-twisted Higgs bundle on $X$ and let $E'$ be a Hecke modification of $E$ determined by $v\in\pp(E_{x}^{\vee})$ at $x\in X$. Assume that the following diagram
$$\xymatrix{0\ar[r]&E'\ar[r]\ar[d]^{\phi'}&E\ar[r]^{v}\ar[d]^{\phi}&\cc_{x}\ar[r]\ar[d]^{\mu}&0\\
0\ar[r]&E'\otimes L\ar[r]&E\otimes L\ar[r]_{v}&\cc_{x}\ar[r]&0}$$
commutes, where each row is a short exact sequence, $\mu:E_{x}/\ker v\to (E\otimes L)_{x}/\ker v$ is induced from the homomorphism $\phi_{x}:E_{x}\to (E\otimes L)_{x}$ and the induced Higgs field $\phi'$ is holomorphic. $(E',\phi')$ is called the \textbf{Hecke modification of }$(E,\phi)$\textbf{ determined by }$v$ at $x$.
\end{definition}

When $L=K_{X}$, G. Wilkin gave a criterion for the existence of a holomorphic Higgs field $\phi'$. It is straightforward to generalize this criterion for twisted Higgs bundles as follows.

\begin{proposition}[Corollary 4.3 of \cite{Wil}]\label{holomorphicity of phi'}
Let $(E,\phi)$ be an $L$-twisted Higgs bundle on $X$ and let $E'$ be a Hecke modification of $E$ determined by $v\in\pp(E_{x}^{\vee})$ at $x\in X$. The following conditions are equivalent
\begin{enumerate}
\item The induced Higgs field $\phi':E'\to E'\otimes L$ is holomorphic.

\item There exists an eigenvalue $\mu$ of $\phi_{x}$ such that $v(\phi(s))=\mu v(s)$ for all sections $s$ of $E$.

\item There exists an eigenvalue $\mu$ of $\phi_{x}$ such that $v$ descends to a well-defined homomorphism $v'\in(\coker(\phi_{x}-\mu\cdot\id))^{\vee}$.
\end{enumerate}
\end{proposition}

He also got the following correspondence between Hecke modifications of Higgs bundles and those of their spectral data in the case that the associated spectral curve is smooth.

\begin{proposition}[Lemma 4.8. of \cite{Wil}]\label{a correspondence between HM on X and HM on the spectral cover}
Let $(E,\phi)$ be a Higgs bundle on $X$. Let $p_{\phi}:X_{\det\phi}\to X$ be the spectral cover and let $\cL_{(E,\phi)}$ be a rank $1$ torsion free sheaf on $X_{\det\phi}$ corresponding to $(E,\phi)$ in the sense of Theorem \ref{BNR correspondence}. Assume that $\det\phi$ has only simple zeros. Then a Hecke modification of $(E,\phi)$ at $x$ corresponds to a Hecke modification of $\cL_{(E,\phi)}$ at $p_{\phi}^{-1}(x)$.
\end{proposition}

The following result is a generalization of Proposition \ref{a correspondence between HM on X and HM on the spectral cover}. We consider the case that $(E,\phi)$ is a twisted Higgs bundle on $X$ and the associated spectral curve $X_{\det\phi}$ is possibly singular and irreducible.

\begin{proposition}\label{a correspondence between HM on X and HM on the normalization of the spectral cover}
Let $(E,\phi)$ be an $L$-twisted Higgs bundle on $X$. Let $p_{\phi}:X_{\det\phi}\to X$ be the spectral cover and let $\cL_{(E,\phi)}$ be a rank $1$ torsion free sheaf on $X_{\det\phi}$ corresponding to $(E,\phi)$ in the sense of Theorem \ref{BNR correspondence}. Let $\tp_{\phi}:\tX_{\det\phi}\to X_{\det\phi}$ be the normalization. Assume that $\cL_{(E,\phi)}$ is locally free and that $\det\phi$ has at least one multiple zero and $\det\phi$ is not a square of a section of $L$ so that $X_{\det\phi}$ is possibly singular and irreducible. Suppose that $X_{\det\phi}$ has $r_{1}$ nodes $x_{1},\cdots,x_{r_{1}}$ of type $A_{m_{i}-1}$ ($i=1,\cdots,r_{1}$) with $m_{i}$ even and $r_{2}$ cusps $y_{1},\cdots,y_{r_{2}}$ of type $A_{m'_{j}-1}$ ($j=1,\cdots,r_{2}$) with $m'_{j}$ odd. Then for each $x$ such that $p_{\phi}^{-1}(x)\not\in\{x_{1},\cdots,x_{r_1},y_{1},\cdots,y_{r_2}\}$, a Hecke modification of $(E,\phi)$ at $x$ corresponds to a Hecke modification of $\tp_{\phi}^{*}\cL_{(E,\phi)}$ over $\tX_{\det\phi}$ at a point in $(p_{\phi}\circ\tp_{\phi})^{-1}(x)$.
\end{proposition}
\begin{proof}
Let $\tp_{\phi}^{-1}(x_{i})=\{x_{i1},x_{i2}\}$ and $\tp_{\phi}^{-1}(y_{j})=\{\ty_{j}\}$. Since $\pr_{1}|_{\tau^{-1}(\Jac(X_{\det\phi}))}$ can be identified with $\tp_{\phi}^{*}:\Jac(X_{\det\phi})\to\Jac(\tX_{\det\phi})$ (See \ref{Parabolic modules on the normalization of the spectral curve}), $\cL_{(E,\phi)}$ corresponds to a parabolic module
$$(\tp_{\phi}^{*}\cL_{(E,\phi)},(U_{1},\cdots,U_{r_{1}},U'_{1},\cdots,U'_{r_{2}})),$$
where $U_{i}$ is an $m_{i}/2$-dimensional subspace of $((\tp_{\phi}^{*}\cL_{(E,\phi)})_{x_{i1}}\oplus(\tp_{\phi}^{*}\cL_{(E,\phi)})_{x_{i2}})^{m_{i}/2}$ and $U'_{j}$ is an $(m'_{j}-1)/2$-dimensional subspace of $((\tp_{\phi}^{*}\cL_{(E,\phi)})_{\ty_{j}})^{m'_{j}-1}$. $\cL_{(E,\phi)}$ fits into the exact sequence
$$0\to\cL_{(E,\phi)}\to\tp_{\phi*}\tp_{\phi}^{*}\cL_{(E,\phi)}\to T\to0,$$
where
$$T=\bigoplus_{i=1}^{r_{1}}\frac{((\tp_{\phi}^{*}\cL_{(E,\phi)})_{x_{i1}}\oplus(\tp_{\phi}^{*}\cL_{(E,\phi)})_{x_{i2}})^{m_{i}/2}}{U_{i}}\oplus\bigoplus_{j=1}^{r_{2}}\frac{((\tp_{\phi}^{*}\cL_{(E,\phi)})_{\ty_{j}})^{m'_{j}-1}}{U'_{j}}.$$
For $\tmu\in\tX_{\det\phi}$ such that $\mu=\tp_{\phi}(\tmu)$ and $x=p_{\phi}(\mu)$, let $\tcL'$ be the Hecke modification of $\tp_{\phi}^{*}\cL_{(E,\phi)}$ determined by $\tv\in\pp((\tp_{\phi}^{*}\cL_{(E,\phi)})_{\tmu}^{\vee})$ at $\tmu$. Then we have the following short exact sequence
$$\xymatrix{0\ar[r]&\tcL'\ar[r]&\tp_{\phi}^{*}\cL_{(E,\phi)}\ar[r]^{\tv}&\cc_{\tmu}\ar[r]&0}$$
and the induced surjective homomorphism
$$q:\tp_{\phi*}\tcL'\to T\to0.$$

Let $\cL'=\ker q$. By the snake lemma, we have the following commutative diagram
\begin{equation}\xymatrix{&0\ar[d]&0\ar[d]&0\ar[d]&\\
0\ar[r]&\cL'\ar[r]\ar[d]_{l}&\tp_{\phi*}\tcL'\ar[r]^{q}\ar[d]&T\ar[r]\ar[d]^{\cong}&0\\
0\ar[r]&\cL_{(E,\phi)}\ar[r]^{\iota}\ar[d]&\tp_{\phi*}\tp_{\phi}^{*}\cL_{(E,\phi)}\ar[r]\ar[d]^{\tp_{\phi*}\tv}&T\ar[r]\ar[d]&0\\
0\ar[r]&\coker l\ar[r]\ar[d]&\tp_{\phi*}\cc_{\tmu}\ar[r]\ar[d]&0\ar[d]&\\
&0&0&0,&}\label{snake}\end{equation}
where both the first and the second row are short exact sequences and the third row is exact.

When $\tmu\not\in\tp_{\phi}^{-1}(\{x_{1},\cdots,x_{r_1},y_{1},\cdots,y_{r_2}\})$, the third row of (\ref{snake}) is
$$0\to\cc_{\mu}\to\cc_{\mu}\to0.$$

Let $v=p_{\phi*}(\tp_{\phi*}\tv\circ\iota):E\to\cc_{x}$ that is surjective. Since the tautological section $\lambda:\cL_{(E,\phi)}\to\cL_{(E,\phi)}\otimes p_{\phi}^{*}L$ induces $\mu:(\cL_{(E,\phi)})_{\mu}\to(\cL_{(E,\phi)})_{\mu}$,
it follows from Theorem \ref{BNR correspondence} that
$$v(\phi(\alpha))=p_{\phi*}(\tp_{\phi*}\tv\circ\iota)(\phi(\alpha))=\tp_{\phi*}\tv\circ\iota(\lambda\alpha)=\mu(\tp_{\phi*}\tv\circ\iota(\alpha))=\mu v(\alpha)$$
for $\alpha\in E(U)=\cL_{(E,\phi)}(p_{\phi}^{-1}(U))$, and then the short exact sequence
$$\xymatrix{0\ar[r]&\cL'\ar[r]^{l}&\cL_{(E,\phi)}\ar[r]^{\tp_{\phi*}\tv\circ\iota}&\cc_{\mu}\ar[r]&0}$$
yields to
$$\xymatrix{0\ar[r]&E'\ar[r]\ar[d]^{\phi'}&E\ar[r]^{v}\ar[d]^{\phi}&\cc_{x}\ar[r]\ar[d]^{\mu}&0\\
0\ar[r]&E'\otimes L\ar[r]&E\otimes L\ar[r]_{v}&\cc_{x}\ar[r]&0}$$
via $p_{\phi*}$. By Proposition \ref{holomorphicity of phi'}, $(E',\phi')$ is the Hecke modification of $(E,\phi)$ determined by $v$ at $x$.

Conversely, let $(E',\phi')$ be the Hecke modification of $(E,\phi)$ determined by $v\in\pp(E_{x}^{\vee})$ at $x\in X$ for $x$ such that $p_{\phi}^{-1}(x)\not\in\{x_{1},\cdots,x_{r_1},y_{1},\cdots,y_{r_2}\}$. By Proposition \ref{holomorphicity of phi'}, $v$ induces a well-defined homomorphism $v'\in(\coker(\phi_{x}-\mu\cdot\id))^{\vee}$. Consider the following diagram
$$\xymatrix{\tp_{\phi}^{*}p_{\phi}^{*}(E\otimes L^{-1})_{\tmu}\ar[rr]^{\quad\tp_{\phi}^{*}((p_{\phi}^{*}\phi-\mu\cdot\id)_{\mu})}&&\tp_{\phi}^{*}p_{\phi}^{*}E_{\tmu}\ar[r]&\tp_{\phi}^{*}\coker(p_{\phi}^{*}\phi-\mu\cdot\id)_{\tmu}\ar[r]&0\\
(E\otimes L^{-1})_{x}\ar[rr]^{\quad(\phi-\mu\cdot\id)_{x}}\ar[u]^{\cong}&&E_{x}\ar[r]\ar[u]^{\cong}&\coker(\phi-\mu\cdot\id)_{x}\ar[r]\ar@{-->}[u]&0}$$
for $\tmu\in\tX_{\det\phi}$ such that $\mu=\tp_{\phi}(\tmu)$ and $x=p_{\phi}(\mu)$. By the universal property of cokernels, there exists a nonzero homomorphism $\coker(\phi-\mu\cdot\id)_{x}\to\tp_{\phi}^{*}\coker(p_{\phi}^{*}\phi-\mu\cdot\id)_{\tmu}$. Since $\dim\coker(\phi-\mu\cdot\id)_{x}=\dim\tp_{\phi}^{*}\coker(p_{\phi}^{*}\phi-\mu\cdot\id)_{\tmu}=1$, this homomorphism is indeed an isomorphism. Since $\cL_{(E,\phi)}=\coker(p_{\phi}^{*}\phi-\mu\cdot\id)$, $v'$ gives a well-defined surjective homomorphism $\tv:\tp_{\phi}^{*}\cL_{(E,\phi)}\to\cc_{\tmu}$ and then we get the Hecke modification of $\tp_{\phi}^{*}\cL_{(E,\phi)}$ determined by $\tv\in\pp((\tp_{\phi}^{*}\cL_{(E,\phi)})_{\tmu}^{\vee})$ at $\tmu$.
\end{proof}

Now we can determine a space parametrizing Hecke modifications associated to a twisted Higgs bundle.

\begin{remark}
Let $(E,\phi)$ be an $L$-twisted Higgs bundle on $X$.
\begin{enumerate}
\item When $X_{\det\phi}$ is smooth, $\pp\cL_{(E,\phi)}^{\vee}$ parametrizes Hecke modifications of $(E,\phi)$.

\item When $X_{\det\phi}$ is singular and irreducible with $r_{1}$ nodes $x_{1},\cdots,x_{r_{1}}$ and $r_{2}$ cusps $y_{1},\cdots,y_{r_{2}}$ as in Proposition \ref{a correspondence between HM on X and HM on the normalization of the spectral cover}, $\pp(\tp_{\phi}^{*}\cL_{(E,\phi)}^{\vee})|_{\tp_{\phi}^{-1}(X_{\det\phi}\setminus\{x_{1},\cdots,x_{r_{1}},y_{1},\cdots,y_{r_{2}}\})}$ parametrizes Hecke modifications of $(E,\phi)$.
\end{enumerate}
\end{remark}

\subsection{$(k,l)$-stability and Hecke modifications of twisted Higgs bundles}

For a nontrivial vector bundle $V$ on $X$ and $k\in\zz$, let $\displaystyle\mu_{k}(V):=\frac{\deg V+k}{\rk(V)}$. Note that any $L$-twisted Higgs bundle on $X$ has rank 2.

\begin{definition}
An $L$-twisted Higgs bundle $(E,\phi)$ on $X$ is \textbf{$(k,l)$-stable} (respectively, \textbf{$(k,l)$-semistable}) if, for every proper $\phi$-invariant line subbundle $F$ of $E$, we have
$$\mu_{k}(F)<\mu_{-l}(E/F)\text{ (respectively, }\mu_{k}(F)\leq\mu_{-l}(E/F)).$$
\end{definition}

\begin{remark}\label{properties of kl stability}
\begin{enumerate}
\item $\mu_{k}(F)<\mu_{-l}(E/F)\text{ (respectively, }\mu_{k}(F)\leq\mu_{-l}(E/F))$ is equivalent to
$$\mu_{k}(F)<\mu_{k-l}(E)\text{ (respectively, }\mu_{k}(F)\leq\mu_{k-l}(E)).$$
\item Recall that $\bM=\bM_{\cO_{X},L}$ is the coarse moduli space of S-equivalence classes of semistable $L$-twisted Higgs bundles with trivial determinant on $X$ (See \ref{twisted Higgs bundles}). For $(E,\phi)\in\bM$, $(E,\phi)$ is $(0,1)$-stable if and only if $(E,\phi)$ is stable.
\end{enumerate}
\end{remark}

\begin{lemma}\label{stability of elm trans}
Fix $x\in X$ and $v\in\pp(E_{x}^{\vee})$. Let $(E,\phi)$ be an $L$-twisted Higgs bundle and let $(E',\phi')$ be a Hecke modification of $(E,\phi)$ determined by $v$ at $x$. Then if $(E,\phi)$ is $(k,l)$-stable, then $(E',\phi')$ is $(k,l-1)$-stable.
\end{lemma}
\begin{proof}
We follow the idea of the proof of Lemma 5.5 in \cite{NR2}. Assume that $(E,\phi)$ is $(k,l)$-stable. Let $F'$ be a $\phi'$-invariant line subbundle of $E'$. Let $F$ be the saturation of the image of the map $F'\to E$. Then $F$ is a $\phi$-invariant line subbundle of $E$ such that $\deg F'\le\deg F$. Since $\deg E=\deg E'+1$, we have
$$\mu_{k}(F')\le\mu_{k}(F)<\mu_{k-l}(E)<\mu_{k-l+1}(E')$$
Hence, $(E',\phi')$ is $(k,l-1)$-stable.
\end{proof}

By Remark \ref{properties of kl stability}-(2) and Lemma \ref{stability of elm trans}, we see that if $(E,\phi)\in\bM$ is stable, then $(E',\phi')$ is stable.

\section{A Hecke cycle associated to a twisted Higgs bundle}\label{Cycles}

In this section, we construct a cycle associated to a stable twisted Higgs bundle and then prove Theorem \ref{Main result1} (Theorem \ref{Main theorem1:map theta}) and Theorem \ref{Main result2} (Theorem \ref{Main theorem2:map Phi}). This cycle is a modified version of the classical Hecke cycle associated to a stable vector bundle (See Definition 5.12 and Theorem 5.13 of \cite{NR2}).

\subsection{The simultaneous resolution of singularities}

It is useful later to introduce the simultaneous resolution of singularities of a family of curves. Consult \cite{DS,Tei} for details.

Let $p:\cC\to S$ be a flat family of reduced curves, where $S$ is any separated scheme.

\begin{definition}[Definition 1.1 of \cite{DS}]
The family $p:\cC\to S$ is \textbf{equigeneric} if
\begin{enumerate}
\item $S$ is reduced,

\item the locus of singular points of fibres is proper over $S$, and

\item $\displaystyle\delta(s)=\sum_{x\in\cC_{s}}\delta_{x}$ is a constant function on $s\in S$, where $\cC_{s}$ is the fiber of $p$ over $s\in S$ and $\delta_{x}=\dim_{\cc}(\nu_{*}\cO_{\tcC_{s}}/\cO_{\cC_{s}})_{x}$ for the normalization map $\nu:\tcC_{s}\to\cC_{s}$.
\end{enumerate}
\end{definition}

\begin{definition}[Definition 1.2 of \cite{DS}]
The family $p:\cC\to S$ is \textbf{equisingular} if there exist
\begin{enumerate}
\item disjoint sections $\s_{1},\cdots,\s_{n}$ of $p$, the union of whose images contains the locus of singular points of fibres, and

\item a proper and birational morphism $\ve:\bar{\cC}\to\cC$, such that
\begin{enumerate}
\item the composition $\bar{p}:=p\circ\ve:\bar{\cC}\to S$ is flat,

\item for every $s\in S$, the induced morphism $\ve_{s}:\bar{\cC}_{s}\to\cC_{s}$ is a resolution of singularities, and

\item for $i=1,\cdots,n$, the induced morphism $\bar{p}:\ve^{-1}(\s_{i}(S))\to S$ is locally (on $\ve^{-1}(\s_{i}(S))$) trivial.
\end{enumerate}
\end{enumerate}
\end{definition}

The following two results will be useful in \ref{A cycle associated to a twisted Higgs bundle and main results}.

\begin{proposition}[Proposition 1.4 of \cite{DS}, II-4.2 of \cite{Tei}]\label{equigeneric to equisingular}
Let $p:\cC\to S$ be an equigeneric family of reduced curves. There exists a dense Zariski-open subset $U\subset S$ such that the restriction $\cC\times_{S}U\to U$ is equisingular.
\end{proposition}

\begin{proposition}[Theorem 1.5 of \cite{DS}, I-Theorem 1.3.2 of \cite{Tei}]\label{The simultaneous resolution of singularities}
Let $p:\cC\to S$ be a flat family of reduced curves, where $\cC$ and $S$ are reduced separated schemes of finite type. If $S$ is normal, then the following two conditions are equivalent:
\begin{enumerate}
\item the family $p:\cC\to S$ is equigeneric;
\item there exists a proper and birational morphism $\ve:\bar{\cC}\to\cC$, such that $\bar{p}=p\circ\ve$ is flat, and for every $s\in S$ , the induced morphism $\bar{\cC}_{s}\to\cC_{s}$ is a resolution of singularities of the fibre $\cC_{s}=p^{-1}(s)$.
\end{enumerate}
In addition, whenever it exists, the simultaneous resolution $\ve$ is necessarily the normalization of $\cC$.
\end{proposition}

\subsection{Moduli space of rational maps from nonsingular curves to a fixed variety}\label{Moduli space of rational maps}

Let $Y$ be a complex variety. In this subsection, we define a stack of rational maps to $Y$.

\begin{definition}
A \textbf{rational map to $Y$} is a triple $(C,U_{C}\subset C,U_{C}\to Y)$, where $C$ is a nonsingular curve, $U_{C}$ is a nonempty open subset of $C$ and $U_{C}\to Y$ is a morphism.
\end{definition}

We define the fibered category $\cR_{Y}$ over the category $(Sch^{sf}/\cc)$ of separated schemes of finite type over $\cc$ whose sections over a scheme $T\in(Sch^{sf}/\cc)$ are triples $(\cC\to T,U_{\cC}\subset \cC,U_{\cC}\to Y\times T)$, where $\cC\to T$ is a flat proper morphism of finite type whose geometric fibers are nonsingular curves, $U_{\cC}$ is a nonempty open subscheme of $\cC$ over $T$ and $U_{\cC}\to Y\times T$ is a morphism over $T$. It is easy to check that $\cR_{Y}$ is a stack over $(Sch^{sf}/\cc)$.

\subsection{A Hecke cycle associated to a twisted Higgs bundle and main results}\label{A cycle associated to a twisted Higgs bundle and main results}

Let $(\cE,\varphi)\to X\times T$ be a family of rank $2$ stable $L$-twisted Higgs bundles on $X$, with fixed determinant $\Lambda$, parametrized by a scheme $T\in(Sch^{sf}/\cc)$. Assume that it gives a family of irreducible spectral covers $p_{\varphi}:X_{\det\varphi}\to X\times T$ and a family of divisors $\cD=\cD_{\det\varphi}=\mathrm{div}(\det\varphi)\subset X\times T$ such that $(p_{\varphi})_{t}=p_{\varphi_{t}}$ and $p_{\varphi_{t}}:X_{\det\varphi_{t}}\to X$ is ramified over $\cD_{t}=D_{\det\varphi_{t}}$ for each $t\in T$. By the correspondence of Theorem \ref{BNR correspondence}, $(\cE,\varphi)$ corresponds to a family $\cL_{(\cE,\varphi)}$ of semistable rank $1$ torsion free sheaves on $X_{\det\varphi}$, parametrized by $T$. Assume that for each $t\in T$, $X_{\det\varphi_{t}}$ has $r_{1}(t)$ nodes of type $A_{m_{i}(t)-1}$ ($i=1,\cdots,r_{1}(t)$) with $m_{i}(t)$ even and $r_{2}(t)$ cusps of type $A_{m'_{j}(t)-1}$ ($j=1,\cdots,r_{2}(t)$) with $m'_{j}(t)$ odd, and $\cL_{(\cE,\varphi)}$ is a family of semistable rank $1$ locally free sheaves.

We observe that $T$ is normal. $\pr_{T}\circ p_{\varphi}:X_{\det\varphi}\to T$ factors uniquely through the normalization $\tT$ of $T$, where $\pr_{T}:X\times T\to T$ is the projection onto $T$. The maps $X_{\det\varphi}\to\tT$ and $\tT\to T$ is denoted by $f$ and $g$ respectively. Since $X_{\det\varphi}$ is a family of irreducible spectral curves, $f^{-1}(\tlt)\cong X_{\det\varphi_{t}}$ for each $\tlt\in\tT$ with $g(\tlt)=t$, and hence $g:\tT\to T$ is an isomorphism.

Since $\pr_{T}\circ p_{\varphi}:X_{\det\varphi}\to T$ is equigeneric (that is, $\displaystyle\sum_{x\in X_{\det\varphi_{t}}}\delta_{x}$ is $2\deg L$ independent of $t\in T$), it follows from Proposition \ref{The simultaneous resolution of singularities} that there exists a proper birational morphism $\tp_{\varphi}:\tX_{\det\varphi}\to X_{\det\varphi}$ such that $\pr_{T}\circ p_{\varphi}\circ\tp_{\varphi}:\tX_{\det\varphi}\to T$ is flat and for each $t\in T$, $\tp_{\varphi}|_{\tX_{\det\varphi_{t}}}$ is identified with the normalization $\tp_{\varphi_{t}}:\tX_{\det\varphi_{t}}\to X_{\det\varphi_{t}}$. Since the genus of $\tX_{\det\varphi_{t}}$ is $2g-1+r_{2}(t)/2$ (Remark 3.2 of \cite{BNR} and Remark 4.3 of \cite{GO}) and $\pr_{T}\circ p_{\varphi}\circ\tp_{\varphi}:\tX_{\det\varphi}\to T$ is flat, $r_{2}(t)$ is independent of $t\in T$. From now on, $r_{2}(t)$ will be denoted by $r_{2}$.

Since $T$ is of finite type and then noetherian, it follows from Proposition \ref{equigeneric to equisingular} that there exist finite disjoint locally closed subsets $T_{1},\cdots,T_{n}$ of $T$ such that $T_{1}$ is dense in $T$, $T_{k+1}$ is dense in $T-T_{k}$ for all $k\ge1$ and $X_{\det\varphi}^{k}=T_{k}\times_{T}X_{\det\varphi}$ is equisingular with $r_{k1}$ nodes $\alpha_{1}^{k},\cdots,\alpha_{r_{k1}}^{k}$ of type $A_{m_{ki}-1}$ ($i=1,\cdots,r_{k1}$) with $m_{ki}$ even and $r_{2}$ cusps $\beta_{1}^{k},\cdots,\beta_{r_{2}}^{k}$ of type $A_{m'_{kj}-1}$ ($j=1,\cdots,r_{2}$) with $m'_{kj}$ odd for each $k$. Let $\tp_{\varphi}^{-1}(\alpha_{i}^{k})=\{\alpha_{i1}^{k},\alpha_{i2}^{k}\}$ and $\tp_{\varphi}^{-1}(\beta_{j}^{k})=\{\tbe_{j}^{k}\}$.

Let $(\cL_{(\cE,\varphi)})_{T_{k}}=(\cL_{(\cE,\varphi)})|_{X_{\det\varphi}^{k}}$, $\tX_{\det\varphi}^{k}=\tp_{\varphi}^{-1}(X_{\det\varphi}^{k})$, $\tp_{\varphi}^{*}(\cL_{(\cE,\varphi)})_{T_{k}}=\tp_{\varphi}^{*}(\cL_{(\cE,\varphi)})|_{\tX_{\det\varphi}^{k}}$ and $$\cD^{k}=\sum_{i=1}^{r_{k1}}\frac{m_{ki}}{2}\alpha_{i1}^{k}+\sum_{i=1}^{r_{k1}}\frac{m_{ki}}{2}\alpha_{i2}^{k}+\sum_{j=1}^{r_2}\frac{m'_{kj}-1}{2}\tbe_{j}^{k}.$$

Relativizing the idea of the proof of Proposition \ref{identification via tau}, we consider $(\tp_{\varphi}^{*}(\cL_{(\cE,\varphi)})_{T_{k}},\cU(\tp_{\varphi}^{*}(\cL_{(\cE,\varphi)})_{T_{k}}))$ as the flat family of parabolic modules on $\tX_{\det\varphi}^{k}$, where $\cU(\tp_{\varphi}^{*}(\cL_{(\cE,\varphi)})_{T_{k}})$ is the subbundle $$(\cL_{(\cE,\varphi)})_{T_{k}}\otimes_{\cO_{X_{\det\varphi}^{k}}}\bigg(\frac{\cO_{X_{\det\varphi}^{k}}}{\tp_{\varphi*}\cO_{\tX_{\det\varphi}^{k}}(-\cD^{k})}\bigg)$$
of
$$\bigg((\cL_{(\cE,\varphi)})_{T_{k}}\otimes_{\cO_{X_{\det\varphi}^{k}}}\bigg(\frac{\cO_{X_{\det\varphi}^{k}}}{\tp_{\varphi*}\cO_{\tX_{\det\varphi}^{k}}(-\cD^{k})}\bigg)\bigg)\otimes_{\cO_{X_{\det\varphi}^{k}}}\tp_{\varphi*}\cO_{\tX_{\det\varphi}^{k}}$$
$$=\tp_{\varphi}^{*}(\cL_{(\cE,\varphi)})_{T_{k}}\otimes_{\cO_{\tX_{\det\varphi}^{k}}}\frac{\cO_{\tX_{\det\varphi}^{k}}}{\cO_{\tX_{\det\varphi}^{k}}(-\cD^{k})}$$
and both $\tp_{\varphi}^{*}(\cL_{(\cE,\varphi)})_{T_{k}}\otimes_{\cO_{\tX_{\det\varphi}^{k}}}\frac{\cO_{\tX_{\det\varphi}^{k}}}{\cO_{\tX_{\det\varphi}^{k}}(-\cD^{k})}$ and $\frac{\tp_{\varphi}^{*}(\cL_{(\cE,\varphi)})_{T_{k}}\otimes_{\cO_{\tX_{\det\varphi}^{k}}}\frac{\cO_{\tX_{\det\varphi}^{k}}}{\cO_{\tX_{\det\varphi}^{k}}(-\cD^{k})}}{\cU(\tp_{\varphi}^{*}(\cL_{(\cE,\varphi)})_{T_{k}})}=(\cL_{(\cE,\varphi)})_{T_{k}}\otimes_{\cO_{X_{\det\varphi}^{k}}}(\tp_{\varphi*}\cO_{\tX_{\det\varphi}^{k}}/\cO_{X_{\det\varphi}^{k}})$ are flat over $T_{k}$.

Let $\pi_{\cL_{(\cE,\varphi)}}:\pp(\tp_{\varphi}^{*}\cL_{(\cE,\varphi)}^{\vee})\to\tX_{\det\varphi}$ be the projection. Then we have the canonical surjective morphism of coherent sheaves
$$s:p_{1}^{*}\tp_{\varphi}^{*}\cL_{(\cE,\varphi)}\to(\pi_{\cL_{(\cE,\varphi)}},\id)_{*}\cO_{\pp(\tp_{\varphi}^{*}\cL_{(\cE,\varphi)}^{\vee})}(1)$$
on $\tX_{\det\varphi}\times_{T}\pp(\tp_{\varphi}^{*}\cL_{(\cE,\varphi)}^{\vee})$, where $\tX_{\det\varphi}\times_{T}\pp(\tp_{\varphi}^{*}\cL_{(\cE,\varphi)}^{\vee})$ is the fibred product of $\pr_{T}\circ p_{\varphi}\circ\tp_{\varphi}:\tX_{\det\varphi}\to T$ and $\pr_{T}\circ p_{\varphi}\circ\tp_{\varphi}\circ\pi_{\cL_{(\cE,\varphi)}}:\pp(\tp_{\varphi}^{*}\cL_{(\cE,\varphi)}^{\vee})\to T$, $(\pi_{\cL_{(\cE,\varphi)}},\id):\pp(\tp_{\varphi}^{*}\cL_{(\cE,\varphi)}^{\vee})\hookrightarrow\tX_{\det\varphi}\times_{T}\pp(\tp_{\varphi}^{*}\cL_{(\cE,\varphi)}^{\vee})$ is a divisor, $p_{1}:\tX_{\det\varphi}\times_{T}\pp(\tp_{\varphi}^{*}\cL_{(\cE,\varphi)}^{\vee})\to\tX_{\det\varphi}$ is the projection onto the first factor and $s$ is the composition of the surjective morphisms
$$p_{1}^{*}\tp_{\varphi}^{*}\cL_{(\cE,\varphi)}\to(\pi_{\cL_{(\cE,\varphi)}},\id)_{*}(\pi_{\cL_{(\cE,\varphi)}},\id)^{*}p_{1}^{*}\tp_{\varphi}^{*}\cL_{(\cE,\varphi)}$$
and
$$(\pi_{\cL_{(\cE,\varphi)}},\id)_{*}\pi_{\cL_{(\cE,\varphi)}}^{*}\tp_{\varphi}^{*}\cL_{(\cE,\varphi)}\to(\pi_{\cL_{(\cE,\varphi)}},\id)_{*}\cO_{\pp(\tp_{\varphi}^{*}\cL_{(\cE,\varphi)}^{\vee})}(1).$$
Let $H(\cL_{(\cE,\varphi)})_{T_{k}}$ be the kernel of the surjective morphism
$$(\tp_{\varphi},\id)_{*}(\ker s)_{T_{k}}\to\frac{(\ker s)_{T_{k}}\otimes_{p_{1}^{*}\cO_{\tX_{\det\varphi}^{k}}}p_{1}^{*}(\cO_{\tX_{\det\varphi}^{k}}/\cO_{\tX_{\det\varphi}^{k}}(-\cD^{k}))}{p_{1}^{*}\cU(\tp_{\varphi}^{*}(\cL_{(\cE,\varphi)})_{T_{k}})\cap((\ker s)_{T_{k}}\otimes_{p_{1}^{*}\cO_{\tX_{\det\varphi}^{k}}}p_{1}^{*}(\cO_{\tX_{\det\varphi}^{k}}/\cO_{\tX_{\det\varphi}^{k}}(-\cD^{k}))),}$$
where $(\ker s)_{T_{k}}=(\ker s)|_{\tX_{\det\varphi}^{k}\times_{T_{k}}\pp(\tp_{\varphi}^{*}\cL_{(\cE,\varphi)}^{\vee})}$. Let $\displaystyle H(\cL_{(\cE,\varphi)})=\bigsqcup_{k=1}^{n}H(\cL_{(\cE,\varphi)})_{T_{k}}$.

\begin{lemma}\label{flatness of H}
$H(\cL_{(\cE,\varphi)})|_{X_{\det\varphi}\times_{T}\pp(\tp_{\varphi}^{*}\cL_{(\cE,\varphi)}^{\vee})|_{\tp_{\varphi}^{-1}(X_{\det\varphi}\setminus\cup_{k=1}^{n}\{\alpha_{i}^{k},\beta_{j}^{k}|1\le i\le r_{k1},1\le j\le r_{2}\})}}$ is flat over $\pp(\tp_{\varphi}^{*}\cL_{(\cE,\varphi)}^{\vee})|_{\tp_{\varphi}^{-1}(X_{\det\varphi}\setminus\cup_{k=1}^{n}\{\alpha_{i}^{k},\beta_{j}^{k}|1\le i\le r_{k1},1\le j\le r_{2}\})}$.
\end{lemma}
\begin{proof}
Over $\cP^{k}:=X_{\det\varphi}^{k}\times_{T_{k}}\pp(\tp_{\varphi}^{*}\cL_{(\cE,\varphi)}^{\vee})|_{\tp_{\varphi}^{-1}(X_{\det\varphi}^{k}\setminus\{\alpha_{i}^{k},\beta_{j}^{k}|1\le i\le r_{k1},1\le j\le r_{2}\})}$, we have the following short exact sequence
$$0\to H(\cL_{(\cE,\varphi)})_{T_{k}}|_{\cP^{k}}\to(\tp_{\varphi},\id)_{*}(\ker s)_{T_{k}}|_{\cP^{k}}\to\cT|_{\cP^{k}}\to0,$$
where
$$\cT=\frac{(\ker s)_{T_{k}}\otimes_{p_{1}^{*}\cO_{\tX_{\det\varphi}^{k}}}p_{1}^{*}(\cO_{\tX_{\det\varphi}^{k}}/\cO_{\tX_{\det\varphi}^{k}}(-\cD^{k}))}{p_{1}^{*}\cU(\tp_{\varphi}^{*}(\cL_{(\cE,\varphi)})_{T_{k}})\cap((\ker s)_{T_{k}}\otimes_{p_{1}^{*}\cO_{\tX_{\det\varphi}^{k}}}p_{1}^{*}(\cO_{\tX_{\det\varphi}^{k}}/\cO_{\tX_{\det\varphi}^{k}}(-\cD^{k}))).}$$
Since $\supp\cD^{k}\cap\tp_{\varphi}^{-1}(X_{\det\varphi}^{k}\setminus\{\alpha_{i}^{k},\beta_{j}^{k}|1\le i\le r_{k1},1\le j\le r_{2}\})=\emptyset$,
$$\cT|_{\cP^{k}}=\frac{p_{1}^{*}\tp_{\varphi}^{*}(\cL_{(\cE,\varphi)})_{T_{k}}\otimes_{p_{1}^{*}\cO_{\tX_{\det\varphi}^{k}}}p_{1}^{*}(\cO_{\tX_{\det\varphi}^{k}}/\cO_{\tX_{\det\varphi}^{k}}(-\cD^{k}))}{p_{1}^{*}\cU(\tp_{\varphi}^{*}(\cL_{(\cE,\varphi)})_{T_{k}})}|_{\cP^{k}}.$$

Since $(\tp_{\varphi},\id)_{*}(\ker s)_{T_{k}}|_{\cP^{k}}$ and
$$\cT|_{\cP^{k}}=\frac{p_{1}^{*}\tp_{\varphi}^{*}(\cL_{(\cE,\varphi)})_{T_{k}}\otimes_{p_{1}^{*}\cO_{\tX_{\det\varphi}^{k}}}p_{1}^{*}(\cO_{\tX_{\det\varphi}^{k}}/\cO_{\tX_{\det\varphi}^{k}}(-\cD^{k}))}{p_{1}^{*}\cU(\tp_{\varphi}^{*}(\cL_{(\cE,\varphi)})_{T_{k}})}|_{\cP^{k}}$$
$$=p_{1}^{*}\bigg((\cL_{(\cE,\varphi)})_{T_{k}}\otimes_{\cO_{X_{\det\varphi}^{k}}}(\tp_{\varphi*}\cO_{\tX_{\det\varphi}^{k}}/\cO_{X_{\det\varphi}^{k}})\bigg)|_{\cP^{k}}$$
are flat over $\pp(\tp_{\varphi}^{*}\cL_{(\cE,\varphi)}^{\vee})|_{\tp_{\varphi}^{-1}(X_{\det\varphi}^{k}\setminus\{\alpha_{i}^{k},\beta_{j}^{k}|1\le i\le r_{k1},1\le j\le r_{2}\})}$, $H(\cL_{(\cE,\varphi)})_{T_{k}}|_{\cP^{k}}$ is also flat over $\pp(\tp_{\varphi}^{*}\cL_{(\cE,\varphi)}^{\vee})|_{\tp_{\varphi}^{-1}(X_{\det\varphi}^{k}\setminus\{\alpha_{i}^{k},\beta_{j}^{k}|1\le i\le r_{k1},1\le j\le r_{2}\})}$.
\end{proof}

Then $H(\cL_{(\cE,\varphi)})|_{X_{\det\varphi}\times_{T}\pp(\tp_{\varphi}^{*}\cL_{(\cE,\varphi)}^{\vee})|_{\tp_{\varphi}^{-1}(X_{\det\varphi}\setminus\cup_{k=1}^{n}\{\alpha_{i}^{k},\beta_{j}^{k}|1\le i\le r_{k1},1\le j\le r_{2}\})}}$ is a flat family of semistable rank 1 torsion free sheaves parametrized by $\pp(\tp_{\varphi}^{*}\cL_{(\cE,\varphi)}^{\vee})|_{\tp_{\varphi}^{-1}(X_{\det\varphi}\setminus\cup_{k=1}^{n}\{\alpha_{i}^{k},\beta_{j}^{k}|1\le i\le r_{k1},1\le j\le r_{2}\})}$.

Let $H(\cE)=(p_{\varphi},\id)_{*}H(\cL_{(\cE,\varphi)})$ and let $H(\varphi)$ be the morphism
$$(p_{\varphi},\id)_{*}H(\cL_{(\cE,\varphi)})\to(p_{\varphi},\id)_{*}H(\cL_{(\cE,\varphi)})\otimes p_{X}^{*}L=(p_{\varphi},\id)_{*}(H(\cL_{(\cE,\varphi)})\otimes(p_{\varphi},\id)^{*}p_{X}^{*}L)$$
induced from the multiplication by the tautological section $\lambda\in p_{\varphi}^{*}\pr_{X}^{*}L$
$$H(\cL_{(\cE,\varphi)})\to H(\cL_{(\cE,\varphi)})\otimes(p_{\varphi},\id)^{*}p_{X}^{*}L=H(\cL_{(\cE,\varphi)})\otimes q_{1}^{*}p_{\varphi}^{*}\pr_{X}^{*}L,$$
where $p_{X}:X\times\pp(\tp_{\varphi}^{*}\cL_{(\cE,\varphi)}^{\vee})\to X$ and $\pr_{X}:X\times T\to X$ are the projections onto $X$ and $q_{1}:X_{\det\varphi}\times_{T}\pp(\tp_{\varphi}^{*}\cL_{(\cE,\varphi)}^{\vee})\to X_{\det\varphi}$ is the projection onto $X_{\det\varphi}$. The following lemma is a simple computation of the determinant of fibers of $H(\cE)$.

\begin{lemma}\label{degree of HM}
$$\det H(\cE)_{\tmu}=\det\cE_{l}\otimes\cO_{X}(-x)$$
for each $\tmu\in\pp(\tp_{\varphi}^{*}\cL_{(\cE,\varphi)}^{\vee})|_{\tp_{\varphi}^{-1}(X_{\det\varphi}\setminus\cup_{k=1}^{n}\{\alpha_{i}^{k},\beta_{j}^{k}|1\le i\le r_{k1},1\le j\le r_{2}\})}$, where $l=\pr_{T}\circ p_{\varphi}\circ\tp_{\varphi}\circ\pi_{\cL_{(\cE,\varphi)}}(\tmu)$ and $x=\pr_{X}\circ p_{\varphi}\circ\tp_{\varphi}\circ\pi_{\cL_{(\cE,\varphi)}}(\tmu)$.
\end{lemma}
\begin{proof}
Let $\mu=\tp_{\varphi}\circ\pi_{\cL_{(\cE,\varphi)}}(\tmu)$. By Remark 3.1 of \cite{BNR} and Proposition 3.10 of \cite{HP}, we have
$$\det\cE_{l}=\det p_{\varphi_{l}*}\cO_{X_{\det\varphi_{l}}}\otimes\Nm_{X_{\det\varphi_{l}}/X}(\cL_{(\cE_{l},\varphi_{l})})=L^{-1}\otimes\Nm_{X_{\det\varphi_{l}}/X}(\cL_{(\cE_{l},\varphi_{l})})$$
and
$$\det H(\cE)_{\tmu}=\det p_{\varphi_{l}*}\cO_{X_{\det\varphi_{l}}}\otimes\Nm_{X_{\det\varphi_{l}}/X}(H(\cL_{(\cE,\varphi)})_{\tmu})=L^{-1}\otimes\Nm_{X_{\det\varphi_{l}}/X}(H(\cL_{(\cE,\varphi)})_{\tmu}).$$
Further, the diagram (\ref{snake}) gives
$$H(\cL_{(\cE,\varphi)})_{\tmu}=\cL_{(\cE_{l},\varphi_{l})}\otimes\cO_{X_{\det\varphi_{l}}}(-\mu).$$
Hence by Proposition 3.1 of \cite{HP} and Lemma 3.4 of \cite{HP},
$$\det H(\cE)_{\tmu}=\det\cE_{l}\otimes\Nm_{X_{\det\varphi_{l}}/X}(\cO_{X_{\det\varphi_{l}}}(-\mu))$$
$$=\det\cE_{l}\otimes\Nm_{\tX_{\det\varphi_{l}}/X}(\cO_{\tX_{\det\varphi_{l}}}(-\tmu))=\det\cE_{l}\otimes\cO_{X}(-x).$$
\end{proof}

By Lemma \ref{stability of elm trans}, Lemma \ref{flatness of H} and Lemma \ref{degree of HM},
$$(H(\cE),H(\varphi))|_{p_{\varphi}(X_{\det\varphi}\setminus\cup_{k=1}^{n}\{\alpha_{i}^{k},\beta_{j}^{k}|1\le i\le r_{k1},1\le j\le r_{2}\})}$$
is a flat family of stable $L$-twisted Higgs bundles on $X$ parametrized by $\pp(\tp_{\varphi}^{*}\cL_{(\cE,\varphi)}^{\vee})|_{\tp_{\varphi}^{-1}(X_{\det\varphi}\setminus\cup_{k=1}^{n}\{\alpha_{i}^{k},\beta_{j}^{k}|1\le i\le r_{k1},1\le j\le r_{2}\})}$ with $\det H(\cE)_{\tmu}=\Lambda\otimes\cO_{X}(-\pr_{X}\circ p_{\varphi}\circ\tp_{\varphi}\circ\pi_{\cL_{(\cE,\varphi)}}(\tmu))$ for each $\tmu\in\pp(\tp_{\varphi}^{*}\cL_{(\cE,\varphi)}^{\vee})|_{\tp_{\varphi}^{-1}(X_{\det\varphi}\setminus\cup_{k=1}^{n}\{\alpha_{i}^{k},\beta_{j}^{k}|1\le i\le r_{k1},1\le j\le r_{2}\})}$. This family provides a rational map
$$\xymatrix{\theta_{(\cE,\varphi)}:\pp(\tp_{\varphi}^{*}\cL_{(\cE,\varphi)}^{\vee})\ar@{-->}[r]&\bM_{X}:=\displaystyle\bigcup_{x\in X}\bM_{\Lambda\otimes\cO_{X}(-x),L}}$$
with a commutative diagram
$$\xymatrix{\pp(\tp_{\varphi}^{*}\cL_{(\cE,\varphi)}^{\vee})|_{\tp_{\varphi}^{-1}(X_{\det\varphi}\setminus\cup_{k=1}^{n}\{\alpha_{i}^{k},\beta_{j}^{k}|1\le i\le r_{k1},1\le j\le r_{2}\})}\ar[rrr]^{\quad\quad\quad\quad\quad\quad\quad\quad(\theta_{(\cE,\varphi)},\pr_{T}\circ p_{\varphi}\circ\tp_{\varphi}\circ\pi_{\cL_{(\cE,\varphi)}})}\ar[d]_{ p_{\varphi}\circ\tp_{\varphi}\circ\pi_{\cL_{(\cE,\varphi)}}}&&&\bM_{X}\times T\ar[d]^{(\alpha,\id)}\\
X\times T\ar[rrr]^{\id}&&&X\times T}$$
where $\alpha:\bM_{X}\to X$ is given by
$$\cO_{X}(-\alpha((E,\phi)))=(\det E)\otimes\Lambda^{-1}\text{ for }(E,\phi)\in\bM_{X}.$$

\begin{theorem}\label{Main theorem1:map theta}
$(\theta_{(\cE,\varphi)},\pr_{T}\circ p_{\varphi}\circ\tp_{\varphi}\circ\pi_{\cL_{(\cE,\varphi)}}):\pp(\tp_{\varphi}^{*}\cL_{(\cE,\varphi)}^{\vee})\dashrightarrow\bM_{X}\times T$ is a closed immersion on $\tp_{\varphi}^{-1}(X_{\det\varphi}\setminus\cup_{k=1}^{n}\{\alpha_{i}^{k},\beta_{j}^{k}|1\le i\le r_{k1},1\le j\le r_{2}\})$. In particular,
$\theta_{(\cE_{t},\varphi_{t})}:\pp(\tp_{\varphi_{t}}^{*}\cL_{(\cE_{t},\varphi_{t})}^{\vee})\dashrightarrow\bM_{X}$ is a closed immersion on $\tp_{\varphi_{t}}^{-1}(X_{\det\varphi_{t}}\setminus\{\textrm{nodes,cusps}\})$ for each $t\in T$.
\end{theorem}
\begin{proof}
Assume that $T$ is a point. An $L$-twisted Higgs bundle $(E,\phi)$ on $X$ corresponds to a rank $1$ torsion free sheaf $\cL_{(E,\phi)}$ on $X_{\det\phi}$ in the sense of Theorem \ref{BNR correspondence}. Assume that $X_{\det\phi}$ has $r_{1}$ nodes $x_{1},\cdots,x_{r_{1}}$ of type $A_{m_{i}-1}$ ($i=1,\cdots,r_{1}$) with $m_{i}$ even and $r_{2}$ cusps $y_{1},\cdots,y_{r_{2}}$ of type $A_{m'_{j}-1}$ ($j=1,\cdots,r_{2}$) with $m'_{j}$ odd. $P$ denotes $\tp_{\phi}^{-1}(X_{\det\phi}\setminus\{x_{1},\cdots,x_{r_1},y_{1},\cdots,y_{r_2}\})$. Note that $\cL_{(E,\phi)}$ corresponds to a parabolic module
$$(\tp_{\phi}^{*}\cL_{(E,\phi)},(U_{1},\cdots,U_{r_{1}},U'_{1},\cdots,U'_{r_{2}})),$$
where $U_{i}$ is an $m_{i}/2$-dimensional subspace of $((\tp_{\phi}^{*}\cL_{(E,\phi)})_{x_{i1}}\oplus(\tp_{\phi}^{*}\cL_{(E,\phi)})_{x_{i2}})^{m_{i}/2}$ and $U'_{j}$ is an $(m'_{j}-1)/2$-dimensional subspace of $((\tp_{\phi}^{*}\cL_{(E,\phi)})_{\ty_{j}})^{m'_{j}-1}$. Let $\tp_{\phi}^{-1}(x_{i})=\{x_{i1},x_{i2}\}$ and $\tp_{\phi}^{-1}(y_{j})=\{\ty_{j}\}$. Let $\pi_{\cL_{(E,\phi)}}:\pp(\tp_{\phi}^{*}\cL_{(E,\phi)}^{\vee})\to\tX_{\det\phi}$ be the projection.

Then we have the following short exact sequences
$$\xymatrix{0\ar[r]&\ker s\ar[r]&p_{1}^{*}\tp_{\phi}^{*}\cL_{(E,\phi)}\ar[r]^{s\quad\quad\quad\quad}&(\pi_{\cL_{(E,\phi)}},\id)_{*}\cO_{\pp(\tp_{\phi}^{*}\cL_{(E,\phi)}^{\vee})}(1)\ar[r]&0.}$$
and
$$\xymatrix{0\ar[r]&H(\cL_{(E,\phi)})\ar[r]&(\tp_{\phi},\id)_{*}\ker s\ar[r]^{q}&&\\
&&\!\!\!\!\!\!\!\!\!\!\!\!\!\!\!\!\!\!\!\!\!\!\!\!\!\!\!\!\!\!\!\!\!\!\displaystyle\bigoplus_{i=1}^{r_{1}}\frac{(\ker s_{x_{i1}}\oplus\ker s_{x_{i2}})^{m_{i}/2}}{V_{i}}\oplus\bigoplus_{j=1}^{r_{2}}\frac{(\ker s_{\ty_{j}})^{m'_{j}-1}}{V'_{j}}\ar[r]&0,&}$$
where $(\pi_{\cL_{(E,\phi)}},\id):\pp(\tp_{\phi}^{*}\cL_{(E,\phi)}^{\vee})\to\tX_{\det\phi}\times\pp(\tp_{\phi}^{*}\cL_{(E,\phi)}^{\vee})$ is a divisor, $p_{1}:\tX_{\det\phi}\times\pp(\tp_{\phi}^{*}\cL_{(E,\phi)}^{\vee})\to\tX_{\det\phi}$ and $p_{2}:\tX_{\det\phi}\times\pp(\tp_{\phi}^{*}\cL_{(E,\phi)}^{\vee})\to\pp(\tp_{\phi}^{*}\cL_{(E,\phi)}^{\vee})$ is the projection onto the first factor, $V_{i}=p_{1}^{*}\tp_{\phi}^{*}U_{i}\cap((\tp_{\phi},\id)_{*}(\ker s))_{x_{i}}$, $V'_{j}=p_{1}^{*}\tp_{\phi}^{*}U'_{j}\cap((\tp_{\phi},\id)_{*}(\ker s))_{y_{j}}$ and $q$ is the composition of the quotient maps
$$(\tp_{\phi},\id)_{*}(\ker s)\to (\tp_{\phi},\id)_{*}(\ker s)\otimes\big(\bigoplus_{i=1}^{r_{1}}\cc(x_{i})\oplus\bigoplus_{j=1}^{r_{2}}\cc(y_{j})\big)$$
and
$$(\tp_{\phi},\id)_{*}(\ker s)\otimes\big(\bigoplus_{i=1}^{r_{1}}\cc(x_{i})\oplus\bigoplus_{j=1}^{r_{2}}\cc(y_{j})\big)\to$$
$$\bigoplus_{i=1}^{r_{1}}\frac{(\ker s_{x_{i1}}\oplus\ker s_{x_{i2}})^{m_{i}/2}}{V_{i}}\oplus\bigoplus_{j=1}^{r_{2}}\frac{(\ker s_{\ty_{j}})^{m'_{j}-1}}{V'_{j}}.$$

Since $\pp(\tp_{\phi}^{*}\cL_{(E,\phi)}^{\vee})$ is a component of the quot scheme $Q:=\Qu_{\tX_{\det\phi}}(\tp_{\phi}^{*}\cL_{(E,\phi)})$, the restricted tangent bundle $T_{Q}|_{P}$ is identified with
$$p_{2*}\HHom(\ker s,(\pi_{\cL_{(E,\phi)}},\id)_{*}\cO_{\pp(\tp_{\phi}^{*}\cL_{(E,\phi)}^{\vee})}(1))|_{P}$$
by Lemma 5.10 of \cite{NR2}, where $p_{2}:\tX_{\det\phi}\times\pp(\tp_{\phi}^{*}\cL_{(E,\phi)}^{\vee})\to\pp(\tp_{\phi}^{*}\cL_{(E,\phi)}^{\vee})$ is the projection onto the second factor.

By Lemma 5.10 of \cite{NR2} again, the infinitesimal deformation map for the family $\ker s|_{\tX_{\det\phi}\times P}$ of coherent sheaves is given, upto sign, by the connecting homomorphism
$$p_{2*}\HHom(\ker s,(\pi_{\cL_{(E,\phi)}},\id)_{*}\cO_{\pp(\tp_{\phi}^{*}\cL_{(E,\phi)}^{\vee})}(1))|_{P}\to R^{1}p_{2*}\cE nd(\ker s|_{\tX_{\det\phi}\times P}).$$
This fits into the following exact sequence
$$0\to p_{2*}\EEnd(\ker s|_{\tX_{\det\phi}\times P})\to p_{2*}\HHom(\ker s|_{\tX_{\det\phi}\times P},p_{1}^{*}\tp_{\phi}^{*}\cL_{(E,\phi)}|_{\tX_{\det\phi}\times P})$$
$$\to p_{2*}\HHom(\ker s,(\pi_{\cL_{(E,\phi)}},\id)_{*}\cO_{\pp(\tp_{\phi}^{*}\cL_{(E,\phi)}^{\vee})}(1))|_{P}\to R^{1}p_{2*}\EEnd(\ker s|_{\tX_{\det\phi}\times P}).$$
Since $\ker s=\tp_{\phi}^{*}H(\cL_{(E,\phi)})$ is a family of line bundles on $\tX_{\det\phi}$, $p_{2*}\EEnd(\ker s|_{\tX_{\det\phi}\times P})\cong\cO_{\pp(\tp_{\phi}^{*}\cL_{(E,\phi)}^{\vee})}|_{P}$. Further, by Lemma 5.6-(ii) of \cite{NR2}, $p_{2*}\HHom(\ker s|_{\tX_{\det\phi}\times P},p_{1}^{*}\tp_{\phi}^{*}\cL_{(E,\phi)}|_{\tX_{\det\phi}\times P})\cong\cO_{\pp(\tp_{\phi}^{*}\cL_{(E,\phi)}^{\vee})}|_{P}$.
Thus the differential of $\theta_{(E,\phi)}$ is injective on $P$. On the other hand, $\theta_{(E,\phi)}$ is itself injective on $P$ by Proposition \ref{a correspondence between HM on X and HM on the normalization of the spectral cover}. Thus we complete the proof.
\end{proof}

Now we are ready to define a cycle associated to a stable twisted Higgs bundle. We define a morphism of stacks $\Phi$ from $\bM^{lf}\cap\bM^{irr}\cap\bM^{s}$ into $\cR_{\bM_X}\times\Pic^{2\deg L}(X)$ consisting of the morphisms $\Phi(T)$ of sections over $T$ given by
$$(\cE,\varphi)\mapsto((\tX_{\det\varphi},\tp_{\varphi}^{-1}(X_{\det\varphi}^{\textrm{ns}})\subset\tX_{\det\varphi},(\theta_{(\cE,\varphi)},\pr_{T}\circ p_{\varphi}\circ\tp_{\varphi}\circ\pi_{\cL_{(\cE,\varphi)}})|_{\tp_{\varphi}^{-1}(X_{\det\varphi}^{\textrm{ns}})}),$$
$$\cO_{X\times T}(\cD_{\det\varphi})),$$
where $\bM=\bM_{\cO_{X},L}$ is the coarse moduli space of S-equivalence classes of semistable $L$-twisted Higgs bundles with trivial determinant on $X$ (See \ref{twisted Higgs bundles}),
$$\bM^{lf}=\{(E,\phi)\in\bM|\cL_{(E,\phi)}\text{ is a rank }1\text{ locally free sheaf on }X_{\det\phi}\},$$
$$\bM^{irr}=\{(E,\phi)\in\bM|X_{\det\phi}\text{ is irreducible}\},\quad\bM^{s}\text{ is the stable locus of }\bM,$$
$$\cR_{\bM_X}\text{ is the stack of rational maps to }\bM_{X}\text{ (See \ref{Moduli space of rational maps})},$$
$$X_{\det\varphi}^{\textrm{ns}}=X_{\det\varphi}\setminus\bigcup_{t\in T}\{\textrm{nodes, cusps of }X_{\det\varphi_{t}}\}$$
and $\cD_{\det\varphi}=\mathrm{div}(\det\varphi)$ is a family of divisors in $X\times T$ (See the first paragraph of \ref{A cycle associated to a twisted Higgs bundle and main results}).

\begin{theorem}\label{Main theorem2:map Phi}
$\Phi$ is injective on objects up to isomorphism.
\end{theorem}
\begin{proof}
It suffices to show that $\Phi(T)$ is injective on objects up to isomorphism for each scheme $T\in(Sch^{sf}/\cc)$. For every pair of $(\cE,\varphi)$ and $(\cF,\psi)$ in $\bM^{lf}\cap\bM^{irr}\cap\bM^{s}(T)$, we have only to show that if isomorphisms $f:X_{\det\varphi}\to X_{\det\psi}$, $\tf:\tX_{\det\varphi}\to\tX_{\det\psi}$ and $\ovf:X\times T\to X\times T$ fit into the following commutative diagram
$$\xymatrix{\pp(\tp_{\varphi}^{*}\cL_{(\cE,\varphi)}^{\vee})\ar[rr]_{\cong}^{\tf}\ar[d]_{\pi_{\cL_{(\cE,\varphi)}}}^{\cong}&&\pp(\tp_{\psi}^{*}\cL_{(\cF,\psi)}^{\vee})\ar[d]_{\cong}^{\pi_{\cL_{(\cF,\psi)}}}\\
\tX_{\det\varphi}\ar[rr]_{\cong}^{\tf}\ar[d]_{\tp_{\varphi}}&&\tX_{\det\psi}\ar[d]^{\tp_{\psi}}\\
X_{\det\varphi}\ar[rr]_{\cong}^{f}\ar[d]_{p_{\varphi}}&&X_{\det\psi}\ar[d]^{p_{\psi}}\\
X\times T\ar[rr]^{\ovf}_{\cong}&&X\times T},$$
$\tf(\tp_{\varphi}^{-1}(X_{\det\varphi}^{\textrm{ns}}))=\tp_{\psi}^{-1}(X_{\det\psi}^{\textrm{ns}})$, $H(\cL_{(\cE,\varphi)})_{\tmu}\cong f^{*}(H(\cL_{(\cF,\psi)})_{\tf(\tmu)})$ for $\tmu\in\pp(\tp_{\varphi}^{*}\cL_{(\cE,\varphi)}^{\vee})|_{\tp_{\varphi}^{-1}(X_{\det\varphi}^{\textrm{ns}})}$ and $\ovf^{*}(\cO_{X\times T}(\cD_{\det\psi}))=\cO_{X\times T}(\cD_{\det\varphi})$, then it gives rise to an equivalence between $(\cE,\varphi)$ and $(\cF,\psi)$ in $\bM^{lf}\cap\bM^{irr}\cap\bM^{s}(T)$.

We claim that $(\cE_{l},\varphi_{l})\cong\ovf^{*}(\cF_{m},\psi_{m})$ for $\tmu\in\pp(\tp_{\varphi}^{*}\cL_{(\cE,\varphi)}^{\vee})|_{\tp_{\varphi}^{-1}(X_{\det\varphi}^{\textrm{ns}})}$, $l=\pr_{T}\circ p_{\varphi}\circ\tp_{\varphi}\circ\pi_{\cL_{(\cE,\varphi)}}(\tmu)$ and $m=\pr_{T}\circ p_{\psi}\circ\tp_{\psi}\circ\pi_{\cL_{(\cF,\psi)}}\circ\tf(\tmu)$. Since
$$H(\cL_{(\cE,\varphi)})_{\tmu}=\cL_{(\cE_{l},\varphi_{l})}\otimes\cO_{X_{\det\varphi_{l}}}(-\tp_{\varphi}\circ\pi_{\cL_{(\cE,\varphi)}}(\tmu))$$
$$=\cL_{(\cE_{l},\varphi_{l})}\otimes f^{*}\cO_{X_{\det\psi_{m}}}(-f\circ\tp_{\varphi}\circ\pi_{\cL_{(\cE,\varphi)}}(\tmu))$$
and
$$H(\cL_{(\cF,\psi)})_{\tf(\tmu)}=\cL_{(\cF_{m},\psi_{m})}\otimes\cO_{X_{\det\psi_{m}}}(-\tp_{\psi}\circ\pi_{\cL_{(\cF,\psi)}}\circ\tf(\tmu))$$
$$=\cL_{(\cF_{m},\psi_{m})}\otimes\cO_{X_{\det\psi_{m}}}(-f\circ\tp_{\varphi}\circ\pi_{\cL_{(\cE,\varphi)}}(\tmu))$$
for $\tmu\in\pp(\tp_{\varphi}^{*}\cL_{(\cE,\varphi)}^{\vee})|_{\tp_{\varphi}^{-1}(X_{\det\varphi}^{\textrm{ns}})}$ and $l=\pr_{T}\circ p_{\varphi}\circ\tp_{\varphi}\circ\pi_{\cL_{(\cE,\varphi)}}(\tmu)$, we have the isomorphism $\cL_{(\cE_{l},\varphi_{l})}\cong f^{*}\cL_{(\cF_{m},\psi_{m})}$.
Taking the pushforward $p_{\varphi_{l}*}$, we get an isomorphism $(\cE_{l},\varphi_{l})\cong\ovf^{*}(\cF_{m},\psi_{m})$.

Hence there exists a line bundle $\cL$ on $T$ such that $(\cE,\varphi)\cong\ovf^{*}(\cF,\varphi)\otimes\pr_{T}^{*}\cL$ (See Lemma 2.5 of \cite{R} and (3.4) of \cite{Tha}).
\end{proof}

Let $\im\Phi$ be the fully faithful subcategory of $\cR_{\bM_X}\times\Pic^{2\deg L}(X)$ whose sections over $T$ are the pairs of the form
$$((\tX_{\det\varphi},\tp_{\varphi}^{-1}(X_{\det\varphi}^{\textrm{ns}})\subset\tX_{\det\varphi},(\theta_{(\cE,\varphi)},\pr_{T}\circ p_{\varphi}\circ\tp_{\varphi}\circ\pi_{\cL_{(\cE,\varphi)}})|_{\tp_{\varphi}^{-1}(X_{\det\varphi}^{\textrm{ns}})}),\cO_{X\times T}(\cD_{\det\varphi})).$$

For each $(E,\phi)\in\bM^{lf}\cap\bM^{irr}\cap\bM^{s}$ corresponding to a geometric point $\Spec\cc\to\bM^{lf}\cap\bM^{irr}\cap\bM^{s}$, the pair of the form
$$((\tX_{\det\phi},\tp_{\phi}^{-1}(X_{\det\phi}\setminus\{\textrm{nodes,cusps}\})\subset\tX_{\det\phi},\theta_{(E,\phi)}|_{\tp_{\phi}^{-1}(X_{\det\phi}\setminus\{\textrm{nodes,cusps}\})}),\cO_{X}(D_{\det\phi})),$$ which is a section of $\im\Phi$ over $\Spec\cc$, is called \textbf{a Hecke cycle associated to $(E,\phi)$}.

\section*{Acknowledgements}
The author thanks Young-Hoon Kiem for helpful comments on an earlier draft. The author also thanks Graeme Wilkin for his inspiring paper \cite{Wil}.

\end{document}